\documentclass[11pt,reqno]{amsart}

\usepackage[scaled=0.92]{helvet}
\usepackage[T1]{fontenc}
\usepackage{textcomp}
\usepackage{amsthm,amsmath,amssymb,amsfonts,amscd,verbatim,latexsym}
\newtheorem{Theorem}{Theorem}[section]
\newtheorem{Lemma}[Theorem]{Lemma}
\newtheorem{Proposition}[Theorem]{Proposition}

\newtheorem{Example}[Theorem]{Example}

\newcommand{\br}[2]{({\mathbf #1} \, {\mathbf #2})}
\newcommand{\Omg}[2]{\Omega_{{\mathbf #1} \, {\mathbf #2}}}
\newcommand{\uu}{\mathbf u}
\newcommand{\uv}{\mathbf v}
\newcommand{\ux}{\mathbf x} 
\newcommand{\uy}{\mathbf y} 
\newcommand{\uz}{\mathbf z} 
\newcommand{\uw}{\mathbf w} 
\newcommand{\ut}{\mathbf t} 
 
\newcommand{\uc}{\mathcal C} 
\newcommand{\sh}{\mathsf h}
\newcommand{\N}{\mathcal N}
\renewcommand{\P}{\mathbf P}
\newcommand{\T}{\mathcal T}
\renewcommand{\O}{\mathcal O}
\newcommand{\C}{\mathcal C}
\newcommand{\Q}{\mathcal Q}
\newcommand{\E}{\mathcal E}
\newcommand{\bE}{\mathbf E}
\newcommand{\Sc}{\mathbb S}
\newcommand{\F}{\mathcal F}
\newcommand{\G}{\mathcal G}
\newcommand{\PHY}{\text{PHY}}
\newcommand{\pc}{\, \circlearrowleft } 
\newcommand{\field}{\Bbbk} 
 
\newcommand{\ra}{\rightarrow}

\newcommand{\lra}{\longrightarrow}

\newcommand{\demo}{\noindent {\sc Proof.}\;}

\begin{document} 
\title[Combinants of a binary pencil]
{On the linear combinants of a binary pencil}
\author[Abdesselam and Chipalkatti]
{Abdelmalek Abdesselam and Jaydeep Chipalkatti}
\maketitle

\bigskip 

\parbox{11.8cm}{ \small 
{\sc Abstract.} 
Let $A,B$ denote binary forms of order $d$, and 
let $\uc_{2r-1} = (A,B)_{2r-1}$  be the sequence of their 
linear combinants for $1 \le r \le \lfloor\frac{d+1}{2}\rfloor$. It is known that 
$\uc_1,\uc_3$ together determine the pencil $\{A + \lambda \, B\}_{\lambda \in \P^1}$, and hence 
indirectly the higher $\C_{2r-1}$. 
In this paper we exhibit explicit formulae for all $r \ge 3$,  which allow us to 
recover $\uc_{2r-1}$ from the knowledge of $\uc_1$ and $\uc_3$. 
The calculations make use of the symbolic method of classical
invariant theory, as well as the quantum theory of angular momentum. 
Our theorem pertains to the plethysm representation $\wedge^2 S_d$ for the group $SL_2$. We give an example for the group $SL_3$ to show that such a result may hold for other 
categories of representations.} 

\bigskip 

\parbox{12cm}{\small
} 

% \tableofcontents 

\section{Introduction} 
This paper is a thematic sequel to~\cite{AC} and~\cite{JC}. 
The problem solved here was originally posed in \cite{JC} (of which a pr{\'e}cis is given below). 
All of the unexplained notation and terminology used in this paper may be found in 
\cite{AC}. The reader is referred to~\cite{Clebsch,Dolgachev,GrYo,Olver} for some 
foundational material in classical invariant theory and the symbolic method. 
The basics of the representation 
theory of $SL_2$ may be found in~\cite[Lecture 11]{FH} and \cite[Chapter~4]{Sturmfels}. 

\subsection{} 
The base field $\field$ will be of characteristic zero. Let $S_d$ denote 
the $(d+1)$-dimensional irreducible representation of the group $SL_2=SL(2,\field)$. 
We identify $S_d$ with the space of (homogeneous) binary $d$-ics in the variables 
$\ux = \{x_1,x_2\}$. 

Given integers $m,n \ge 0$ and $0 \le q \le \min(m,n)$, there is an 
$SL_2$-equivariant split surjection (see~\cite[\S 1.5]{AC})
\begin{equation} 
\pi_q: S_m \otimes S_n \lra S_{m+n-2q}. \label{pi.q} \end{equation} 
Given binary forms $F \in S_m$ and $G \in S_n$, the image $\pi_q(F \otimes G)$ is classically 
referred to as the $q$-th transvectant of $F$ and $G$, denoted by $(F,G)_q$. We have an 
explicit formula 
\[ (F,G)_q = \frac{(m-q)! \, (n-q)!}{m! \, n!} \, \sum\limits_{i=0}^q \, (-1)^i \, \binom{q}{i} \, 
\frac{\partial^q F}{\partial x_1^{q-i} \, \partial x_2^i} \, 
\frac{\partial^q G}{\partial x_1^i \, \partial x_2^{q-i}} \, ;   \] 
however, it is seldom directly useful. For later use, let 
\begin{equation} 
\imath_q: S_{m+n-2q} \lra S_m \otimes S_n  \label{i.q} \end{equation} 
denote the canonical inclusion, so that $\pi_q \circ \imath_q$ is the identity map on $S_{m+n-2q}$. 
\subsection{} 
Now let $A,B \in S_d$ denote two linearly independent forms. 
There is an isomorphism of $SL_2$-representations 
\begin{equation} 
\wedge^2 S_d = \bigoplus\limits_{r=1}^{\lfloor\frac{d+1}{2}\rfloor} \,  S_{2d-4r+2}, 
\label{wedge2.Sd} \end{equation} 
with projection morphisms $p_r: \wedge^2 S_d \lra S_{2d-4r+2}$. 
The image $p_r(A \wedge B)$ equals the transvectant $(A,B)_{2r-1}$. For any scalars 
$\alpha,\beta,\gamma,\delta$, we have an invariance property 
\[ (\alpha \, A + \beta \, B,\gamma \, A + \delta \, B)_{2r-1} = 
(\alpha \, \delta - \beta \, \gamma) \, (A,B)_{2r-1}. \] 
Hence, up to a scalar, the forms $\uc_{2r-1} = (A,B)_{2r-1}$ 
depend only on the subspace $\Pi_{A,B} = \text{Span} \, \{A,B\}$. 
In classical terminology (see~\cite[\S 250]{GrYo}), the 
$\{\uc_{2r-1}\}$ are linear combinants of the pencil $\{A + \lambda \, B\}_{\lambda \in \P^1}$. 

Decomposition~(\ref{wedge2.Sd}) implies that 
the pencil is completely determined by the sequence of forms 
\[ \uc_1,\; \uc_3,\dots, \uc_{\lfloor \frac{d+1}{2} \rfloor},   \] 
but rather more can be said. An arbitrary form $F \in S_d$ belongs 
$\Pi_{A,B}$, if and only if the Wronskian 
\[ {\mathbf W} = \left| \begin{array}{ccc} 
A_{x_1^2} & A_{x_1x_2} & A_{x_2^2} \\ 
B_{x_1^2} & B_{x_1x_2} & B_{x_2^2} \\ 
F_{x_1^2} & F_{x_1x_2} & F_{x_2^2}  \end{array} \right| = 0. \] 
After some manipulation, this condition can be rewritten as 
\[ (\uc_1,F)_2 + \frac{d-2}{4d-6} \, F \, \uc_3 = 0. \] 
It follows that $\uc_1,\uc_3$ determine $\Pi_{A,B}$, and hence they indirectly determine 
all the subsequent combinants $\uc_5,\uc_7,\uc_9$ etc. It is natural to 
enquire whether there exists a concrete formula for $\uc_{2r-1}$ in terms of 
$\uc_1,\uc_3$. This problem was solved in~\cite[\S 5]{JC} 
for $\uc_5$ and $\uc_7$ using some ad-hoc 
calculations; here we will give an inductive solution which applies to all $r \ge 3$. 

\begin{Example} \rm 
Assume $d=7$. We have an identity 
\begin{equation}
 \uc_1 \, \uc_5 = - \frac{21}{2} \, (\uc_1,\uc_1)_4 +\frac{84}{11} \, 
(\uc_1,\uc_3)_2 + \frac{735}{484} \, \uc_3^2, 
\label{syzygy.d7.r3} \end{equation}
which expresses $\uc_5$ in terms of $\uc_1,\uc_3$. 
Similarly, the identity 
\begin{equation} \begin{aligned} 
\uc_1 \, \uc_7 = & -28 \, (\uc_1,\uc_1)_6 - \frac{210}{11} \, (\uc_1,\uc_3)_4+ 
8 \, (\uc_1,\uc_5)_2 \\ 
& + \, \frac{1960}{121} \, (\uc_3,\uc_3)_2+\frac{35}{11} \, \uc_3 \, \uc_5, 
\end{aligned} 
\label{syzygy.d7.r5} \end{equation}
indirectly expresses $\uc_7$ in terms of $\uc_1,\uc_3$. 
We will show that such formulae always exist for all $d$ and 
$3 \le r \le \lfloor\frac{d+1}{2}\rfloor$. 
\end{Example} 

After completing our results, we discovered that a few such 
calculations had been done by Shenton~\cite[p.~257ff]{Shenton}. 

\section{Quadratic syzygies} 
Define a (quadratic) syzygy of weight $2 r$ to be an identity 
\begin{equation} \sum \, 
\alpha_{i,j} \, (\uc_{2i-1},\uc_{2j-1})_{2(r-i-j+1)} =0, 
\qquad (\alpha_{i,j} \in {\mathbf Q}) \label{quad.syzygy} \end{equation}
assumed to hold for all $d$-ics $A,B$. The sum is quantified over all 
pairs $(i,j)$ such that 
\begin{equation} 1 \le i \le j \le r, \quad i + j \le r+1. 
\label{range.ij} \end{equation} 
For instance, (\ref{syzygy.d7.r3}) and (\ref{syzygy.d7.r5}) are 
syzygies of weight $6$ and $8$ respectively. 
Notice that the only term in~(\ref{quad.syzygy}) involving 
$\uc_{2r-1}$ corresponds to $(i,j)=(1,r)$. Now our main result is the following: 
\begin{Theorem} \sl 
For every $3 \le r \le \lfloor\frac{d+1}{2}\rfloor$,
there exists a quadratic syzygy of 
weight $2 r$ such that $\alpha_{1,r} \neq 0$. 
\label{main.theorem} \end{Theorem} 
We will, in fact, produce an explicit formula for the $\alpha_{i,j}$. 
Given this, one can rewrite (\ref{quad.syzygy}) as 
\[ \uc_{2r-1} = - \frac{1}{\uc_1} \, 
\sum\limits \; \frac{\alpha_{i,j}}{\alpha_{1,r}} \; 
(\uc_{2i-1},\uc_{2j-1})_{2(r-i-j+1)} =0, \] 
which recovers $\uc_{2r-1}$ from $\uc_1,\dots,\uc_{2r-3}$. Notice that 
$\uc_1$ is (up to a scalar) the Jacobian of $A,B$; in particular it 
is nonzero if $\{A,B\}$ are linearly independent. 

By a classical theorem of Gordan, 
the algebra of all combinants of a pencil is finitely generated. However, a specific 
set of generators is known in only a few cases
(see~\cite{Gordan,Meulien,Newstead,Wall}). Our main theorem 
is not directly comparable to these results, since we allow not only polynomial, but 
also rational transvectant expressions in the combinants.

\subsection{} 
In outline, the proof of Theorem~\ref{main.theorem} proceeds as follows. The following proposition (proved in 
\cite[\S5]{JC}) reinterprets a syzygy as an $SL_2$-equivariant morphism. 

\begin{Proposition} \sl \label{prop.wedge4Sd}
The vector space of syzygies of weight $2r$ is isomorphic 
to $\text{Hom}_{SL_2}(S_{4(d-r)},\wedge^4 S_d)$. 
\end{Proposition} 
In \S\ref{zetadef} we will construct a {\sl specific} morphism 
\[ \zeta: S_{4(d-r)} \lra \wedge^4 S_d, \] 
and then calculate the corresponding syzygy coefficients. In fact this calculation will be done 
twice: first by classical symbolic methods, and secondly by recasting the coefficient as a 9-j symbol 
in the sense of the quantum theory of angular momentum. 

It would be of interest to know whether an analogue of Theorem~\ref{main.theorem} 
holds for other categories of representations. 
In \S\ref{example.SL3} we give such an example for the group $SL_3$. 
\subsection{} We informally sketch the idea behind Proposition~\ref{prop.wedge4Sd}. 
Consider the Pl{\"u}cker imbedding 
\[ G(2,S_d) \hookrightarrow \P(\wedge^2 S_d), \] 
with image $X$ and ideal sheaf ${\mathcal I}_X$. The 
short exact sequence of $SL_2$-representations 
\[ 0 \ra H^0({\mathcal I}_X(2)) \ra H^0(\O_\P(2)) 
\ra H^0(\O_X(2)) \ra 0, \] 
can be naturally identified with 
\[ 0 \ra \wedge^4 S_d \stackrel{\imath}{\ra} S_2(\wedge^2 S_d) 
\stackrel{q}{\ra} \Sc_{(2,2)}(S_d) \ra 0. \] 
Here $\Sc_{(2,2)}$ denotes the Schur functor associated to the partition $(2,2)$ 
(see~\cite[Lecture 6]{FH}). The coefficients of each $\uc_{2i-1}$ can be seen as 
homogeneous co{\"o}rdinates on $\P(\wedge^2 S_d)$, hence an 
expression 
\[ \E = \sum \, \alpha_{i,j} \,(\uc_{2i-1},\uc_{2j-1})_{2(r-i-j+1)} \] 
corresponds to the function 
\[ S_{4(d-r)} \stackrel{\phi_\E}{\lra} H^0(\O_\P(2)), \quad 
F \lra (F,\E)_{4(d-r)}. \] 
Now $\E$ is a syzygy iff this function is identically zero 
on $X$, i.e., iff $q \circ \phi_\E=0$. 
This is equivalent to the condition that $\phi_\E$ factor through $\ker q$. Conversely, a nonzero map 
$S_{4(d-r)} \stackrel{\phi}{\lra} \wedge^4 S_d$ defines an irreducible subrepresentation 
of $H^0({\mathcal I}_X(2))$, which translates into a quadratic syzygy $\E_\phi$. \qed 

\subsection{} 
This interpretation allows to read off the individual coefficients in a syzygy. 
Let $\E=0$ denote a quadratic syzygy of weight $2r$, and fix a pair of integers $(i,j)$ 
satisfying 
\[ 1 \le i, j \le r, \quad i + j \le r+1. \]  
(Notice that we have not imposed the condition $i \le j$.) Consider the sequence of morphisms 
\begin{equation} \begin{aligned} 
S_{4(d-r)} & \stackrel{\phi_\E}{\lra} \wedge^4 S_d \stackrel{\imath}{\lra}
S_2(\wedge^2 S_d) \stackrel{\beta_1}{\lra} \wedge^2 S_d \otimes \wedge^2 S_d \\ 
& \stackrel{\beta_2}{\lra} S_{2d-4i+2} \otimes S_{2d-4j+2} 
\stackrel{\beta_3}{\lra} S_{4(d-r)}. 
\end{aligned} \label{beta.seq} \end{equation} 
Here $\beta_1$ is the natural inclusion 
map $v \cdot w \lra \frac{1}{2}(v \otimes w + w \otimes v)$, 
$\beta_2$ is the tensor product of projections $p_{2i-1} \otimes p_{2j-1}$, and $\beta_3$ is the 
transvectant map $\pi_{2(r-i-j+1)}$. By Schur's lemma, the composite endomorphism 
\[ \beta_3 \circ \beta_2 \circ \beta_1 \circ \imath \circ \phi_\E: 
S_{4(d-r)} \lra S_{4(d-r)} \] 
must be the multiplication by a constant, say $\theta_{i,j}$. Then, up to a global constant, 
\begin{equation} 
\E = \sum\limits \; \theta_{i,j} \, (\uc_{2i-1},\uc_{2j-1})_{2(r-i-j+1)}. 
\label{expr.E.theta} \end{equation}

\subsection{} 
In this section we will describe the $\beta_i$ using the 
classical symbolic calculus.  
Our notation follows~\cite{AC} and \cite{GrYo}; in particular, $\ux=(x_1,x_2),\uy=(y_1,y_2)$ etc.~denote 
binary variables, and 
\[ \br{x}{y} = x_1 \, y_2 - y_1 \, x_2, \quad 
\Omg{x}{y} = 
\frac{\partial^2}{\partial x_1 \, \partial y_2}-
\frac{\partial^2}{\partial y_1 \, \partial x_2}. \] 
Define 
\[ \sh(m,n;q) = \frac{(m+n-2q+1)!}{(m+n-q+1)! \, q \, !} \, . \] 
The rationale for introducing this factor is explained in~\cite[\S 1.6]{AC}. 

We will realise $S_2(\wedge^2 S_d)$ as the space of quadrihomogeneous 
forms $Q(\ux,\uy,\uz,\uw)$ of order $d$ in each variable, satisfying the 
conditions 
\[ Q(\ux,\uy,\uz,\uw) = - \, Q(\uy,\ux,\uz,\uw) = - \, Q(\ux,\uy,\uw,\uz) = 
Q(\uz,\uw,\ux,\uy). \] 
Inside this space, the image of $\imath$ is identified with the set of 
alternating forms, i.e., those $Q$ for which 
\[ Q(\ux,\uy,\uz,\uw) = \text{sign}(\sigma) \, 
Q(\ux^\sigma,\uy^\sigma,\uz^\sigma,\uw^\sigma), \] 
for every permutation $\sigma$ of the four letters. 

Now realise $S_{2d-4i+2} \otimes S_{2d-4j+2}$ as the space
of bihomogeneous forms of respective orders $(2d-4i+2, 2d-4j+2)$ in $\uu,\uv$. 
Then $\beta_2 \circ \beta_1$ maps $Q$ to 
\[ \sh(d,d;2i-1) \, \sh(d,d;2j-1) \, 
\left[ \, \Omg{x}{y}^{2i-1} \, \Omg{z}{w}^{2j-1} \ Q \, \right] \] 
followed by the substitutions $\{ \ux,\uy \ra \uu \}, \{ \uz,\uw \ra \uv \}$. 
Notice that, given the two pairs of operations 
\[ \Omg{x}{y}, \{\ux,\uy \ra \uu\}, \qquad \Omg{z}{w}, \{\uz,\uw \ra \uv\}, \] 
any operation from the first pair commutes from any operation from the second. 
Finally realise $S_{4(d-r)}$ as order $4(d-r)$ forms in $\ut$, then $\beta_3$ maps $R(\uu,\uv)$ to 
\[ \sh(2d-4i+2,2d-4j+2;2r-2i-2j+2) \, 
[ \, \Omg{u}{v}^{2r-2i-2j+2} \, R(\uu,\uv) \, ], \] 
followed by the substitutions $\{\uu,\uv \ra \ut\}$.

\subsection{} \label{zetadef}
Now define $\zeta: S_{4(d-r)} \lra S_2(\wedge^2 S_d)$ to be the morphism 
which sends $f_\ut^{4(d-r)}$ to the form 
\[ \begin{aligned} \F = \; 
  & \br{x}{y} \, \br{z}{w}^{2r-1} \, f_\ux^{d-1} \, f_\uy^{d-1} \, f_\uz^{d-2r+1} \, f_\uw^{d-2r+1} \\ 
- \; & \br{x}{z} \, \br{y}{w}^{2r-1} \, f_\ux^{d-1} \, f_\uz^{d-1} \, f_\uy^{d-2r+1} \, f_\uw^{d-2r+1} \\
  & \br{x}{w} \, \br{y}{z}^{2r-1} \, f_\ux^{d-1} \, f_\uw^{d-1} \, f_\uy^{d-2r+1} \, f_\uz^{d-2r+1} \\ 
- \; & \br{y}{w} \, \br{x}{z}^{2r-1} \, f_\uy^{d-1} \, f_\uw^{d-1} \, f_\ux^{d-2r+1} \, f_\uz^{d-2r+1} \\ 
  & \br{z}{w} \, \br{x}{y}^{2r-1} \, f_\uz^{d-1} \, f_\uw^{d-1} \, f_\ux^{d-2r+1} \, f_\uy^{d-2r+1} \\ 
-  \; & \br{z}{y} \, \br{x}{w}^{2r-1} \, f_\uz^{d-1} \, f_\uy^{d-1} \, f_\ux^{d-2r+1} \, f_\uw^{d-2r+1} \, .
\end{aligned} \] 
By construction, $\F$ is alternating in all four variables; hence $\zeta$ factors through $\wedge^4 S_d$. 
The rationale behind this choice of $\zeta$ will be explained in~\S\ref{pos.HS}. 
\subsection{The first calculation} Let us write (using the obvious notation) 
\[ \F = \T(\ux\uy,\uz \uw) - \, \T(\ux\uz,\uy\uw) \, + \dots  - \, \T(\uz\uy, \ux\uw).\] 
We should like to gauge the effect of the morphism 
$\beta_3 \circ \beta_2 \circ \beta_1$ on each summand in $\F$. 
The next two lemmata allow us to `cancel' an $\Omg{x}{y}$ against an
$(\ux \, \uy)$. 

\begin{Lemma} \sl
Let $\G$ denotes an arbitrary bihomogeneous form of orders $p,q$ in $\ux,\uy$ respectively.
\begin{enumerate} 
\item[(a)]For every $m\ge 1$, 
\[
\Omg{x}{y} \, (\ux \, \uy)^m \, \G =
m \, (p+q+m+1) \, (\ux\, \uy)^{m-1} \, \G+ (\ux\, \uy)^m \, \Omg{x}{y} \, \G .
\]
\item[(b)] For every $\ell \ge 1$, 
\[ \Omg{x}{y}^\ell \, (\ux \, \uy) \, \G =
\ell \, (p+q-\ell+3) \, \Omg{x}{y}^{\ell-1} \, \G+ (\ux\, \uy) \, \Omg{x}{y}^\ell \, \G . \] 
\end{enumerate} \label{cancel.lemma} \end{Lemma} 
\demo By straightforward differentiation, 
\[ \begin{aligned} 
{} & \Omg{x}{y}  \, (\ux \, \uy) \, \G = 2 \, \G + 
 (x_1 \, \frac{\partial \G}{\partial x_1} +  x_2 \, \frac{\partial \G}{\partial x_2}) +  
(y_1 \, \frac{\partial \G}{\partial y_1} +  y_2 \, \frac{\partial \G}{\partial y_2}) + \\ 
& (x_1 \, y_2 - x_2 \, y_1) \, 
(\frac{\partial^2 \G}{\partial x_1 \, \partial y_2} - 
\frac{\partial^2 \G}{\partial x_2 \, \partial y_1}) \\ 
= \;  & (p+q+2) \, \G + (\ux \, \uy) \, \Omg{x}{y} \, \G. 
\end{aligned} \] 
Now part (a) follows by an easy induction on $m$, and (b) by one on $\ell$. 
\qed 

\begin{Lemma} \sl With $\G$ as above, and $\ell,m \ge 0$, 
\[ [\, \Omg{x}{y}^\ell \, (\ux \, \uy)^m \, \G \, ]_{\ux,\uy \ra \uu} = 
\begin{cases} 
\mu(p,q;\ell,m) \; 
[\, \Omg{x}{y}^{\ell-m}  \, \G \, ]_{\ux,\uy \ra \uu}  & \text{if $\ell \ge m$}, \\ 
0 & \text{otherwise,} \end{cases} \] 
where 
\begin{equation} \mu(p,q;\ell,m) = \frac{\ell \, !}{(\ell-m)!} 
\frac{(p+q-\ell+2m+1)!}{(p+q-\ell+m+1)!} \, . 
\label{eqn.c3} \end{equation} 
\label{lemma.Ev2} \end{Lemma} 
\demo 
Using part (a) of the previous lemma for the connecting step, one shows by induction on $\ell$, that 
\[
\, \Omg{x}{y}^\ell \, (\ux \, \uy)^m \, \G  \equiv
\begin{cases} 
\mu(p,q;\ell,m) \; 
\Omg{x}{y}^{\ell-m}  \, \G   
& \text{if $0\le \ell-m \le\min(p,q)$}, \\ 
0 & \text{otherwise,} \end{cases} \] 
where $\equiv$ stands for congruence modulo $(\ux\, \uy)$. The result follows, because 
terms involving $(\ux \, \uy)$ vanish after the substitution $\ux,\uy\ra\uu$. \qed 

\medskip 

As a consequence, the term $\T(\uz \uw,\ux \uy)$ is annihilated by 
the operation $\Omg{x}{y}^{2i-1}$ followed by $\{\ux,\uy \ra \uu\}$, 
unless $i=r$ (and hence necessarily $j=1$). In the latter case, 
\[ \begin{aligned} {} & \beta_2 \circ \beta_1 \, (\T(\uz\uw,\ux\uy)) = 
\\ & \sh(d,d,2r-1)\, \sh(d,d,1) \, 
[ \, \Omg{x}{y}^{2r-1} \, \Omg{z}{w} \circ \T(\uz\uw,\ux\uy) \, ]_{\{\ux,\uy \ra \uu\}, \{\uz,\uw \ra \uv\}}, 
\end{aligned} \] 
evaluates to 
 \begin{equation} f_\uu^{2d-4r+2} \, f_\uv^{2d-2}  \label{fuv} \end{equation}
because
\[ \begin{aligned} 
{} & \sh(d,d,2r-1) \; \sh(d,d,1) \; \mu(d-2r+1,d-2r+1,2r-1,2r-1) \, \times \\ 
& \mu(d-1,d-1,1,1) = 1. 
\end{aligned} \] 
Then $\beta_3$ carries (\ref{fuv}) into $f_\ut^{4(d-r)}$. By the same argument, 
$\T(\ux\uy,\uz\uw)$ goes to $f_\ut^{4(d-r)}$ if $(i,j)=(1,r)$, and zero otherwise. 

This disposes of two of the summands in $\F$; the rest of them will need more work. 
As an interlude, we will consider a  preparatory example which illustrates the operation 
of $\Omega_{\ux\uy}$ on a symbolic product involving $\ux,\uy$ (cf.~\cite[\S 3.2.5]{Glenn}). 
\begin{Example} \rm 
Let 
\[ E = \br{x}{z}^7 \, \br{y}{w}^2 \, f_\ux\, g_\ux^4 \, f_\uy^5. \] 
First we follow the calculation of $\Omega_{\ux\uy}  \, E$. 
The idea, in brief, is to pair an $\ux$-factor with a $\uy$-factor and contract them against each other. 
The following diagram shows all the types of $\ux$ and $\uy$ factors in $E$, and the possible 
pairings between them. 

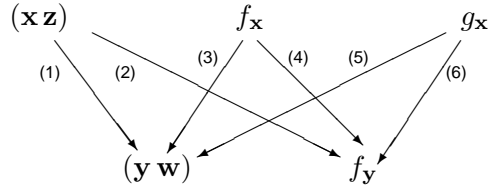
\setlength{\unitlength}{1mm} 
\begin{picture}(100,30)
\put(30,20){$\br{x}{z}$} 
\put(60,20){$f_\ux$} 
\put(90,20){$g_\ux$} 
\put(45,0){$\br{y}{w}$} 
\put(75,0){$f_\uy$} 
%%%%%%%%%%%%%%%
\put(36,18){\vector(3,-4){11}}
\put(41,19){\vector(2,-1){33}}
\put(61,18){\vector(-2,-3){10}}
\put(88,19){\vector(-2,-1){33}}
\put(90,18){\vector(-2,-3){11}}
\put(63,18){\vector(1,-1){14}}
%%%%%%%%%%%%%%%
\put(34,13){\tiny (1)}
\put(44,13){\tiny (2)}
\put(55,15){\tiny (3)}
\put(67,15){\tiny (4)}
\put(75,15){\tiny (5)}
\put(88,13){\tiny (6)}
\end{picture} 

\medskip 

The equality $\Omega_{\ux \uy} \, [ \br{x}{u} \br{y}{v} ] = \br{u}{v}$ gives our basic 
rule: contracting $\br{x}{u}$ against $\br{y}{v}$ gives $\br{u}{v}$. 
For instance, contraction along the arrow (1) 
gives $\br{z}{w}$. Introducing a phantom letter $\tilde{\mathbf f} = (-f_2,f_1)$, we can 
write $f_\ux = (\ux \, \tilde{\mathbf f})$, and hence contraction along (3) 
gives $(\tilde{\mathbf f} \, \uw) = - \, f_\uw$. Contraction along (4) 
gives $(f \, f)=0$. Now $\Omega_{\ux \uy}  \, E$ is a sum of terms 
(quantified over all choices of contractions), where in each term the contracted factors are replaced by 
their result. Thus, $\Omega_{\ux \uy} \, E = $
\[ \underbrace{14 \, \br{z}{w} \, \br{x}{z}^6 \, \br{y}{w} \, f_\ux \, g_\ux^4 \, f_\uy^5}_{\text{from arrow (1)}}
+ \underbrace{35 \, f_\uz \, \br{x}{z}^6 \, \br{y}{w}^2 \, f_\ux\, g_\ux^4 \, f_\uy^4}_{\text{from arrow (2)}}  \, + 
\dots \text{etc.} \] 

To calculate $\Omega_{\ux\uy}^2  \, E$ we must sum over all possible $2$-step sequences 
of contractions, taking account of available multiplicities.  For instance, 
the sequences of arrows 
\[ (1)(4), \quad (2)(2), \quad (3)(5) \] 
are allowed, but (3)(3) is not since there is only one $f_\ux$ available. This gives 
$\Omega_{\ux \uy}^2  \, E =$ 
\[ \underbrace{840 \, f_\uz^2 \, \br{x}{z}^5 \, \br{y}{w}^2 \, f_\ux \, 
g_\ux^4 \, f_\uy^3}_{\text{from (2)(2)}} - 
\underbrace{120 \, (g \, f) \, g_\uw \, \br{x}{z}^7 \, \br{y}{w} \, f_\ux\, g_\ux^2 \, f_\uy^4}_{\text{from (5)(6)}} 
+ \dots \text{etc.} \] 
If we treat the seven $\br{x}{z}$ factors as notionally distinct, a sequence of two from them can 
be chosen in $7!/5!$ ways, and similarly for $f_\uy^5$. This gives the first coefficient as 
$\frac{7!}{5!} \, \frac{5!}{3!} = 840$. Similarly, the second coefficient is $\frac{4!}{2!} \times 2 \times 5$. 
Notice that the sequence (6)(5) will give an {\sl additional} term identical to the 
one coming from (5)(6). 
\label{example.evalOmega} \end{Example} 

\subsection{} We will now follow the evaluation of $\beta_3 \circ \beta_2 \circ \beta_1 \circ \T(\ux\uw,\uy\uz)$. 

As a first step we have to remove $(2i-1)$ factors each of type $\ux,\uy$ from $\T(\ux\uw,\uy\uz)$. 
The available factors are respectively 
\[ \br{x}{w} \, f_\ux^{d-1} \quad \text{and} \quad 
\br{y}{z}^{2r-1} \, f_\uy^{d-2r+1}. \] There are three choices: 
\begin{equation} \begin{array}{crcl} 
\text{(I)} & \br{x}{w} \, f_\ux^{2i-2} & \text{and} & \br{y}{z}^{2i-2} \, f_\uy \\ 
\text{(II)} & \br{x}{w} \, f_\ux^{2i-2} & \text{and} & \br{y}{z}^{2i-1}, \\ 
\text{(III)} & f_\ux^{2i-1}& \text{and} & \br{y}{z}^{2i-1}. 
\end{array} \label{choices.xy} \end{equation}
The possibilities are limited by the following constraint: 
since $f_\uy$ can only be 
paired with $\br{x}{w}$, no more than one copy of $f_\uy$ can be chosen; and hence 
at least $2i-2$ copies of $\br{y}{z}$ must be chosen. 

After contraction and the substitution $\{\ux,\uy \ra \uu\}$, choice (I) leads to the expression 
\begin{equation} 
c_I  \, \br{u}{z}^{2r-2i+1} \, f_\uu^{2d-2r-2i+1} \, f_\uz^{d-2r+2i-1} \, f_\uw^d. 
\label{Omegaxy.I} \end{equation} 
Here (and subsequently) $c_I,c_{I'}$ etc.~stand for some rational constants 
which will be determined later. Now we must remove 
$(2j-1)$ factors each of type $\uz,\uw$ from~(\ref{Omegaxy.I}). The choice is forced, namely 
\[ \begin{array}{crcl} 
\text{(I')} & \br{u}{z}^{2j-1} & \text{and} & f_\uw^{2j-1}. 
\end{array} \] 
After contraction and $\{\uz,\uw \ra \uv\}$, we get an expression 
\begin{equation} 
 - \, c_I \, c_{I'} \, \br{u}{v}^{2r-2i-2j+2} \, f_\uu^{2d-2r-2i+2j} \, 
f_\uv^{2d-2r+2i-2j}. \label{expr.uvI}\end{equation}
(The negative sign arises, because contracting $\br{u}{z}$ against $f_\uw$ gives $-f_\uu$.) 
Now $\beta_3$ will convert~(\ref{expr.uvI}) into 
\begin{equation} - \, c_I \, c_{I'} \, f_\ut^{4(d-r)} 
\label{expr.choiceI} \end{equation}
as a consequence of Lemma \ref{lemma.Ev2}.
\subsection{} 
Choice (II) in (\ref{choices.xy}) leads to the expression 
\[ - \, c_{II} \, \br{z}{w} \, 
\underbrace{\br{u}{z}^{2r-2i} \, f_\uu^{2d-2r-2i+2} \, 
f_\uw^{d-1} \, f_\uz^{d-2r+2i-1}}_{\G}, \] 
on which we have to operate on by $\Omg{z}{w}^{2j-1}$. 
Using part (b) of Lemma~\ref{cancel.lemma}, 
\[
\Omg{z}{w}^{2j-1} \, (\uz\, \uw) \, \G=
(2j-1)(2d-2j+2) \, \Omg{z}{w}^{2j-2} \, \G +(\uz\, \uw) \, \Omg{z}{w}^{2j-1} \, \G. \]
After the substitution $\{\uz,\uw \ra \uv \}$, the second term goes away. 
In evaluating $\Omg{z}{w}^{2j-2} \,  \G$, we have a forced choice 
\[ \begin{array}{clcr} 
\text{(II')} & \br{u}{z}^{2j-2} & \text{and} & f_\uw^{2j-2}, 
\end{array} \] 
leading to 
\begin{equation} - \, (2d-2j+2)(2j-1) \,
c_{II} \, c_{II'} \,  f_\ut^{4(d-r)}. 
\label{expr.choiceII} \end{equation}
\subsection{} 
Choice (III) (which is only possible if $2i\le d$), leads to
\[ -c_{III} \, (\uu\uw) \, (\uu\uz)^{2r-2i} \, f_\uu^{2d-2r-2i+1} \, f_\uw^{d-1} \, f_\uz^{d-2r+2i}.
\]
When applying $\Omg{z}{w}^{2j-1}$, it further bifurcates into the two choices:
\[ \begin{array}{clcr} 
\text{(III')} & \br{u}{z}^{2j-2} \, f_\uz & \text{and} & \br{u}{w} \, f_\uw^{2j-2} \\ 
\text{(III'')} & \br{u}{z}^{2j-1}  & \text{and} & f_\uw^{2j-1},  \\ 
\end{array} \] 
which are dealt with similarly. In fact (III') can arise only if 
\begin{equation} 
2j \le d, \quad \text{and} \quad r \ge i+j. 
\label{condition.CIIIp} \end{equation}
Altogether we arrive at the expression 
\[ \begin{aligned} {} & \beta_3 \circ \beta_2 \circ \beta_1 \circ \T(\ux\uw,\uy\uz)  
= \sh(d,d;2i-1) \, \sh(d,d;2j-1) \, \times \\ 
& (-c_I \, c_{I'} - \, (2d-2j+2)(2j-1) \, c_{II} \, c_{II'} 
- c_{III} \, c_{III'} + c_{III} \, c_{III''}) \, f_\ut^{4(d-r)}.
\end{aligned} \] 
Using the recipe of Example~\ref{example.evalOmega}, we get the constants 
\[ \begin{aligned} 
c_I & = (2i-1) \, (d-2r+1) \, \frac{(d-1)!}{(d-2i+1)!} \frac{(2r-1)!}{(2r-2i+1)!}, \\ 
c_{I'} & = \frac{(2r-2i+1)! \, d!}{(2r-2i-2j+2)! \, (d-2j+1)!}, \\ 
c_{II} & = \frac{(2i-1) \, (d-1)! \, (2r-1)! }{(d-2i+1)! \, (2r-2i)!}, \\ 
c_{II'} & = \frac{(2r-2i)! \, (d-1)!}{(2r-2i-2j+2)! \, (d-2j+1)!}, \\ 
c_{III} & = \frac{(d-1)! \, (2r-1)!}{(d-2i)! \, (2r-2i)!}, \\ 
c_{III'} & = \frac{(2j-1) \, (d-2r+2i) \, (2r-2i)! \, (d-1)!}{(2r-2i-2j+2)! \, (d-2j+1)!}, \\ 
c_{III''} & = \frac{(2r-2i)! \, (d-1)!}{(2r-2i-2j+1)! \, (d-2j)!}. 
\end{aligned} \] 
If $2i\le d$ fails, then $c_{III}$ is zero by definition.
Likewise, if the conditions in~(\ref{condition.CIIIp}) are not satisfied, then 
$c_{III''}$ is understood to be zero. Recall that the prevailing hypotheses are 
\[
3\le r\le \frac{d+1}{2}, \quad 1\le i,j\le r,
\quad \text{and} \quad
i+j\le r+1.
\]
Therefore, any of extra conditions $2i\le d$, $2j\le d$, and $i+j\le r$
can only fail if respectively $d-2i+1$, $d-2j+1$, or $r-i-j+1$ vanish. 
Hence, the following expressions for $c_{III}$ and $c_{III''}$ hold unconditionally: 
\begin{eqnarray} 
c_{III} & = & \frac{(d-2i+1)\, (d-1)! \, (2r-1)!}{(d-2i+1)! \,
  (2r-2i)!}, 
\label{improvedIII}\\  
c_{III''} & = & \frac{(d-2j+1)\,(2r-2i-2j+2)\,(2r-2i)! \, (d-1)!}{(2r-2i-2j+2)! \, (d-2j+1)!}. \label{improvedIIIpp}
\end{eqnarray}
Due to the symmetry in the situation, the rest of the terms 
\[ - \, \T(\ux\uz,\uy\uw), \quad - \, \T(\uy\uw,\ux\uz), 
\quad - \, \T(\uz\uy,\ux\uw) \] 
give identical evaluations. After some simplification, we arrive at the following formula: 

\subsection{} 
Define $\delta_{i,j}$ to be $1$ if $i=j$, and $0$ otherwise. Let 
\[ \begin{aligned} 
\N_1 = & \, (2 d \, i+2d \, j-d \, r-2i^2-2j^2-2d+3i+3j-2) \, \times \\ 
& d! \, (d-1)! \, (2r-1)! \, (2d-4i+3)! \, (2d-4j+3)!, \\ 
\N_2 = & \, (2i-1)! \, (2j-1)! \, (d-2i+1)! \, (d-2j+1)! \, \times \\ 
& (2d-2i+2)! \, (2d-2j+2)! \, (2r-2i-2j+2)!, 
\end{aligned} \] 
then
\begin{equation} 
\theta_{i,j} = (\, \delta_{i,1} \, \delta_{j,r} + 
\delta_{i,r} \, \delta_{j,1} - 8 \; \frac{\N_1}{\N_2} \, ). 
\label{formula.alpha.ij} \end{equation}
Evidently $\theta_{i,j} = \theta_{j,i}$. Therefore, in expression (\ref{expr.E.theta}) one can combine the 
terms $(i,j)$ and $(j,i)$. Let $\epsilon_{i,j} = 2$ if $i \neq j$, and $1$ if $i = j$. Now let 
$\alpha_{i,j} = \epsilon_{i,j} \, \theta_{i,j}$. We have finally arrived at the required syzygy 
\begin{equation} 
\E_{\zeta}: \sum\limits_{(i,j)} \, \alpha_{i,j} \, (\uc_{2i-1},\uc_{2j-1})_{2(r-i-j+1)} =0,
\end{equation}
where the sum is quantified over all pairs $(i,j)$ such that 
\[ 1 \le i \le j \le r, \quad i + j \le r+1. \] 
The reader may check that for $d=7,r=3$, the syzygy becomes 
\[ 10 \, (\uc_1,\uc_1)_4 - \frac{80}{11} \, (\uc_1,\uc_3)_2  
- \frac{175}{121} \, \uc_3^2 + \frac{20}{21} \, \uc_1 \, \uc_5 = 0, 
\] 
which is the same as (\ref{syzygy.d7.r3}). We have (successfully) tested 
formula~(\ref{formula.alpha.ij}) in {\sc Maple} on several examples. 

\subsection{Second calculation} 
In fact, formula~(\ref{formula.alpha.ij}) was first arrived at by a different path, namely 
by interpreting $\theta_{i,j}$ as (in essence) a 9-j symbol in the sense of the quantum theory of 
angular momentum (see~\cite[\S7]{AC}). 

We pick up the thread at the beginning of \S\ref{zetadef}. The trajectory 
$f_\ut^{4(d-r)} \lra \T(\ux \uw, \uy \uz)$ followed by $\beta_3 \circ \beta_2 \circ \beta_1$ is described 
by the sequence of morphisms 
\[  \begin{aligned} 
{} & S_{4(d-r)} \lra S_{2d-2} \otimes S_{2d-4r+2} \lra 
 (S_d \otimes S_d) \otimes (S_d \otimes S_d) \lra \\ 
& (S_d \otimes S_d) \otimes (S_d \otimes S_d) \lra 
S_{2d-4i+2} \otimes S_{2d-4j+2} \lra  S_{4(d-r)}. \end{aligned} \] 
Here the first two maps are natural injections, the last two are natural projections, and the one in 
the middle is the shuffling map 
\[ (v_1 \otimes v_2) \otimes (v_3 \otimes v_4) \lra (v_1 \otimes v_4) \otimes (v_2 \otimes v_3). \] 
By Schur's lemma, the total composite must a multiple of the identity map 
$\text{Id}_{S_{4(d-r)}}$. Up to an 
easily calculated factor (see~\cite[\S7.9]{AC}), this multiple is the 9-j symbol 
\[ B = \left\{ \begin{array}{ccc}
\frac{d}{2} & \frac{d}{2}  & d-2i+1 \\ & & \\
\frac{d}{2} & \frac{d}{2}  & d-2j+1 \\ & & \\ d-1 & d-2r+1 & 2d-2r
\end{array} \right\} \ . \]
Now interchange rows $1,2$ of $B$, then interchange rows $1,3$ of the new array, and finally 
interchange columns $2,3$. This gives an equivalent array 
\begin{equation}
B' = \left\{ \begin{array}{ccc}
d-1 & 2d-2r & d-2r+1\\
& & \\
\frac{d}{2} & d-2i+1 & \frac{d}{2}\\
& & \\
\frac{d}{2} & d-2j+1 & \frac{d}{2}\\
\end{array}\right\}.  \label{fortriplesum} \end{equation}
Finally apply the Ali\v{s}auskas-Jucys triple sum formula (see~\cite[\S 7.10]{AC}) to 
$B'$. In the notation used there, the set $\Lambda$ of triples of
indices which appear in the sum is contained in
\[ \left\{ \, (d-2r+1,2j-1,0), \; (d-2r+1,2j-2,0), \; (d-2r,2j-1,0) \, \right\},  \] 
which reduces the sum to at most three easily manageable terms.
The triple $(d-2r+1,2j-1,0)$ appears in the sum unless
$i=r=\frac{d+1}{2}$.
The triple $(d-2r+1,2j-2,0)$ always appears.
Finally $(d-2r,2j-1,0)$ appears unless $r=\frac{d+1}{2}$.
One can remove the case discussion using the same trick which led to
the unconditional formulae (\ref{improvedIII}) and (\ref{improvedIIIpp}). 
After a little simplification, 
once again we get formula (\ref{formula.alpha.ij}). \qed 

\section{Positivity} 
\subsection{} \label{pos.HS} 
The next proposition will conclude the proof of Theorem~\ref{main.theorem}. 
\begin{Proposition}\label{nonzero.alpha}
\sl The coefficient $\alpha_{1,r}$ is nonzero, and in fact strictly positive. 
\end{Proposition} 
\demo 
We will extensively use the material in~\cite[\S 7]{AC}. If 
$u: {\mathcal E}_1 \lra {\mathcal E}_2$ denotes a linear map between
Hilbert spaces, then 
$u^*: {\mathcal E}_2 \lra {\mathcal E}_1$ denotes its adjoint. Recall that the Hilbert-Schmidt norm 
of $u$ is defined to be 
\[ || u||_{\text{HS}} = \sqrt{\, \text{trace} \, (u^*  \circ u)}. \] 
For a composite ${\mathcal E}_1 \stackrel{u}{\lra} {\mathcal E}_2 \stackrel{v}{\lra} 
{\mathcal E}_3$, we have $(v \circ u)^* = u^* \circ v^*$. 

In the notation of~\cite[\S 7]{AC}, we write ${\mathcal H}_{\frac{m}{2}}$
for $S_m$, which carries a natural structure of a finite dimensional Hilbert space. 
We will view $\zeta$ as a map from ${\mathcal H}_{2(d-r)}$
to $({\mathcal H}_{\frac{d}{2}})^{\otimes 4}$ via the natural inclusion
$\wedge^4 \, {\mathcal H}_{\frac{d}{2}}
\hookrightarrow ({\mathcal H}_{\frac{d}{2}})^{\otimes 4}$. 
Similarly, we view $\beta_1$ as originating from
$({\mathcal H}_{\frac{d}{2}})^{\otimes 4}$ via the natural surjection
\[
\begin{aligned}
({\mathcal H}_{\frac{d}{2}})^{\otimes 4} &
\lra S_2(\wedge^2\, {\mathcal H}_{\frac{d}{2}})\\
z_1\otimes z_2 \otimes z_3\otimes z_4 & \lra
(z_1\wedge z_2)\cdot(z_3\wedge z_4).
\end{aligned}
\]

Henceforth, throughout the proof, 
the symbol $\pc$ will stand for some strictly positive constant which need not be specified. 
Recall that we have defined maps $\pi^\PHY,\imath^\PHY$ such that 
\[ \pi_{\frac{m}{2},\frac{n}{2},\frac{1}{2}(m+n-2q)}^\PHY = \pc \pi_q, \quad 
\imath_{\frac{m}{2},\frac{n}{2},\frac{1}{2}(m+n-2q)}^\PHY = \pc \imath_q, \] 
in the notation of (\ref{pi.q}) and (\ref{i.q}); moreover $\pi^\PHY = (\imath^\PHY)^*$. 
We will show that, 
\begin{equation} 
 \pc \alpha_{1,r} = ||\zeta||_{\text{HS}}^2. \label{alpha1r.zeta} \end{equation} 
First, observe that $\alpha_{1,1} = 2 \, (r-2) (2r-1) \neq 0$, hence the map $\zeta$ is 
not identically zero (if the reader was not already so persuaded).  If $\left( a_{s,t}\right)$ 
denotes the matrix representing $\zeta$ with respect to some
orthonormal
bases, then 
\[ \text{trace} \, (\zeta^* \circ \zeta) = \sum\limits_{s,t} \; |a_{s,t}|^2 > 0,  \] 
hence it only remains to show~(\ref{alpha1r.zeta}) to complete the proof of the proposition. 

Now specialise to $i=1,j=r$, and let 
$\psi = \beta_3 \circ \beta_2 \circ \beta_1 \circ \zeta$. 
By definition, $\alpha_{1,r}=2 \, \theta_{1,r}$, where 
$\psi = \theta_{1,r} \, \text{Id}_{S_{4(d-r)}}$ and hence 
\[ \theta_{1,r} = \frac{\text{trace}(\psi)}{4(d-r)+1}. \] 
Notice that, up to a positive multiplicative constant, the map $f_\ut^{4(d-r)} \lra \T(\ux \uy,\uz \uw)$ is the sequence 
\[ {\mathcal H}_{2(d-r)} \lra {\mathcal H}_{d-1} \otimes
{\mathcal H}_{d-2r+1} \lra
({\mathcal H}_{\frac{d}{2}} \otimes {\mathcal H}_{\frac{d}{2}}) \otimes 
({\mathcal H}_{\frac{d}{2}} \otimes {\mathcal H}_{\frac{d}{2}} ), \] 
where the first map is $\imath_{j_{12}j_{34}J}^{\PHY}$, and the second is
$\imath_{j_1 j_2 j_{12}}^{\PHY}\otimes \imath_{j_3 j_4 j_{34}}^{\PHY}$,
with
\[
\begin{aligned} 
{} & j_{12} = d-1, \quad j_{34} = d-2r+1, \quad J = 2(d-r), \quad \text{and} \\ 
{} & j_1 = j_2 = j_3 = j_4 = \frac{d}{2}. 
\end{aligned} \]
If we compose this with the alternation map 
\[ \begin{aligned} 
{\mathcal A}: ({\mathcal H}_{\frac{d}{2}})^{\otimes 4} & \lra 
({\mathcal H}_{\frac{d}{2}})^{\otimes 4} \\ 
z_1 \otimes z_2 \otimes z_3 \otimes z_4 & \lra 
\frac{1}{4!}\sum\limits_{\sigma \in {\mathfrak S}_4} \; 
\text{sign}(\sigma) \, z_{\sigma(1)} \otimes z_{\sigma(2)} \otimes z_{\sigma(3)} \otimes z_{\sigma(4)}, 
\end{aligned} \] 
the net effect (up to a constant) is $\zeta: f_\ut^{4(d-r)} \lra {\mathcal F}$. 
In other words, 
\[ \zeta = \pc \, {\mathcal A} \circ \left( 
\imath^\PHY_{j_1 j_2 j_{12}} \otimes \imath^\PHY_{j_3 j_4 j_{34}} \right) 
\circ \imath^\PHY_{j_{12} j_{34} J}. \] 
Now observe that 
$\beta_3 \circ \beta_2 \circ \beta_1$ is (up to a constant) 
the sequence of maps:
\[
({\mathcal H}_{\frac{d}{2}} \otimes {\mathcal H}_{\frac{d}{2}}) \otimes 
({\mathcal H}_{\frac{d}{2}} \otimes {\mathcal H}_{\frac{d}{2}})
\lra {\mathcal H}_{d-1} \otimes {\mathcal H}_{d-2r+1} \lra 
{\mathcal H}_{2(d-r)}, \]
where the first map is 
$\pi_{j_1 j_2 j_{12}}^{\PHY}\otimes \pi_{j_3 j_4 j_{34}}^{\PHY}$, and the second is 
$\pi^\PHY_{j_{12} j_{34} J}$. 
Since the 
maps $\imath^\PHY$ and $\pi^\PHY$ (with identical subscripts) are mutually 
adjoint, and ${\mathcal A}$ is a self-adjoint idempotent,
\[ \psi = \pc \, \zeta^* \circ \zeta, \] 
and the claim follows. \qed 

\medskip 

Indeed, it was this argument which led us to the correct guess for
$\F$.
One strategy to ensure 
that $\alpha_{1,r}$ does not vanish is to make it appear as the
Hilbert-Schmidt norm of a nonzero operator. 
This prompted us to take the adjoint of $\beta_3\circ\beta_2\circ\beta_1$, which determines the first term in $\F$ and 
hence all the rest. 

\subsection{} 
The result of Proposition \ref{nonzero.alpha} amounts to the inequality
\begin{equation}
4 \, (dr-2r^2+3r-1) \times\frac{(d-1)!\,(2d-4r+3)!}{(d-2r+1)!\,(2d-2r+2)!}  <1 
\label{elemineq} \end{equation} 
in the range $r \ge 3, d \ge 2r-1$. 

We include an elementary proof of this inequality. Let $\Gamma(r,d)$ denote the left-hand side of 
(\ref{elemineq}). First, 
\[ \Gamma(r,2r-1) = \frac{2}{r} < 1.   \] 
Let us write 
\[ \frac{\Gamma(r,d+1)}{\Gamma(r,d)} = \frac{N}{D}, \] 
where 
\[ \begin{aligned} 
N & =   d \, (d \, r + 4 \, r - 2 \, r^2-1) \, (2d-4r+5), \\ 
D & = (d \, r - 2 \, r^2 + 3 \, r -1) \, (2d-2r+3) \, (d-r+2). 
\end{aligned} \] 
Now observe that 
\[ D-N = (r-1) \, (r-2) \, (2r-1) \, (d-2r+3) > 0, \] 
hence $\Gamma(r,d+1) < \Gamma(r,d)$. 
This completes the proof. \qed 

\medskip 

Unfortunately this proof gives no insight into why the inequality should be true. It seems 
especially fortuitous that $D-N$  should admit such a tidy factorisation. For reasons already 
stated, we prefer the earlier argument. 

Note that the Hilbert-Schmidt idea also guided the construction
of the closed form syzygy in~\cite[\S 2.14]{AC}. It can be used to
provide an alternate proof of Lemma~2.3 therein.
In \cite{AC0} and \cite[Proposition 5]{SB} one may find similar instances, 
where the nonvanishing of an algebraic expression 
produced by a tensorial construction is the key ingredient in a geometric result.

\section{A ternary example}  \label{example.SL3} 
Our main theorem leads to the analogous problem for $SL_N$-representations. 
To wit, let $V$ denote an $N$-dimensional vector space and write $\Sc_\lambda$ for the 
Schur module $\Sc_\lambda \, V$ (see~\cite[Lecture 6]{FH}). 
Assume that we are given a plethysm decomposition\footnote{To the best of our knowledge, 
no explicit formula for the multiplicities $M_\nu$ is known for an arbitrary $\lambda$. 
See~\cite[Ch.~I.8]{MacDonald} for some special cases.} of Schur modules
\begin{equation}  \wedge^2 \, \Sc_\lambda \simeq 
\bigoplus\limits_\nu \; (\Sc_\nu \otimes \field^{M_\nu}). \label{decom.SLN} \end{equation} 
Let $\mathcal C = \{\uc_\nu^{(i)}: \quad 1 \le i \le M_\nu \}$
denote the associated linear combinants of 
a pencil of tensors in $\Sc_\lambda$. 
It is a natural problem to find a subcollection of ${\mathcal C}$ 
which determines the rest of them. We will now exhibit such an example in the 
ternary case. The symbolic formalism used below is explained in~\cite[\S 4]{AC}. 
\subsection{} Assume 
$N=3$ and $\lambda = (3,1)$. We have a decomposition 
\[ \wedge^2 \, \Sc_{(3,1)} \simeq 
\underbrace{\Sc_{(5)} \oplus \Sc_{(5,3)} \oplus \Sc_{(4,1)} \oplus 
\Sc_{(3,2)} \oplus \Sc_{(1,1)}}_{\bE}, \] 
with projection morphisms 
$f_\lambda: \wedge^2 \, \Sc_{(3,1)} \ra \Sc_\lambda$. 
Let 
\[ A = (a \, b \, \uu) \, a_\ux^2, \quad 
B = (c \, d \, \uu) \, c_\ux^2, \] 
denote two `generic' forms in $\Sc_{(3,1)}$, 
and write $\uc_\lambda = f_\lambda(A \wedge B)$. Then we have symbolic formulae 
\[ \begin{aligned}
\uc_{(5)} & = (a \, b \, d) \, a_\ux^2 \, c_\ux^3, \\ 
\uc_{(5,3)} & = (a \, b \, \uu) \,  (a \, c \, \uu) \, (a \, d \, \uu) \, c_\ux^2, \\ 
\uc_{(4,1)} & = (a \,b \,d) (a \, c \, \uu) \, a_\ux \, c_\ux^2 
-5 \, (a \,b \, c) (a \, d \, \uu) \, a_\ux \, c_\ux^2, \\ 
\uc_{(3,2)} & = (a \,b \,d) \, (a \,c \, \uu)^2 \, c_\ux + 
(a \,b \,c) \, (a \,c \,\uu) \, (a \, d \, \uu) \, c_\ux, \\ 
\uc_{(1,1)} & = (a \,b \,c) \, (a \,c \,d) \, (a \,c \, \uu). 
\end{aligned} \] 
There is an exact sequence of $SL_3$-representations 
\[ 0 \ra \underbrace{\wedge^4 \, \Sc_{(3,1)}}_\Q \ra \Sc_{(2)}(\wedge^2 \, \Sc_{(3,1)}) 
\ra \Sc_{(2,2)}(\Sc_{(3,1)}) \ra 0, \] and, as in the binary case, 
the irreducible subrepresentations of $\Q$ correspond to the 
quadratic syzygies between the $\uc_\lambda$. 

\begin{Proposition} \sl 
Either of the combinants $\uc_{(3,2)}$ and $\uc_{(1,1)}$ can be recovered 
from the set $\{ \uc_{(5)},\uc_{(5,3)},\uc_{(4,1)} \}$. 
\end{Proposition} 
The result follows from an explicit calculation involving plethysms and projection maps. 
Taking our cue from the binary case, we look for subrepresentations corresponding to 
$(5,0) + (3,2) = (8,2)$. 
Decomposing\footnote{The full decompositions are very lengthy, and it seems needless to list them here. 
All plethysm decomposition throughout this 
example were calculated using the `SF' (Symmetric Functions) package for 
{\sc Maple} written by John Stembridge.} 
$\Q$ and $\Sc_2(\bE)$ into irreducible summands, 
we found that they respectively contain $2$ and $7$ copies of $\Sc_{(8,2)}$. The 
latter come from tensor products of the summands in $\bE$ taken two at a time;  e.g., 
the morphism 
\[ \Sc_{(5)} \otimes \Sc_{(5,3)} \lra \Sc_{(8,2)} \] 
is given by the formula 
\[ a_\ux^5 \otimes (c \, d \, \uu)^3 \, c_\ux^2 \ra 
(a \, c \, d) \, (a \, d \, \uu)^2 \, a_\ux^2 \, c_\ux^4. \] 
Let us write $\langle \uc_{(5)}, \uc_{(5,3)}\rangle$ for the 
image of $\uc_{(5)} \otimes \uc_{(5,3)}$ via this morphism. Once all the seven maps have been 
written down symbolically, it only remains to solve a system of linear equations to find the 
two-dimensional space of syzygies; this 
was done in {\sc Maple}.  One conveniently chosen syzygy is the following: 
\begin{equation}  \begin{aligned} 
\uc_{(5)} \, \uc_{(3,2)} = \, 
& \frac{1}{7680} \, \langle \uc_{(5)}, \uc_{(5)} \rangle + 
\frac{1}{92160} \, \langle \uc_{(5)}, \uc_{(5,3)} \rangle -
\frac{1}{5760} \, \langle \uc_{(5)}, \uc_{(4,1)} \rangle \\ 
- & \frac{1}{204800} \, \langle \uc_{(5,3)}, \uc_{(5,3)} \rangle -
\frac{1}{51200} \, \langle \uc_{(5,3)}, \uc_{(4,1)} \rangle. 
\end{aligned} \label{C5.C32} \end{equation} 
This gives a formula for $\uc_{(3,2)}$ in terms of 
$\uc_{(5)}, \uc_{(5,3)}, \uc_{(4,1)}$. 

\smallskip 

There are respectively $3$ and $9$ copies of $\Sc_{(6,1)}$ in 
$\Q$ and $\Sc_2(\bE)$, and the corresponding syzygies are found similarly. 
The following syzygy 
\begin{equation} \begin{aligned} \uc_{(5)} \, \uc_{(1,1)} = \, 
& \frac{1}{25920} \, \langle \uc_{(5)}, \uc_{(5,3)} \rangle - 
\frac{1}{1728} \, \langle \uc_{(5)}, \uc_{(4,1)} \rangle -
\frac{1}{864} \, \langle \uc_{(5)}, \uc_{(3,2)} \rangle \\ 
- & \frac{11}{69120} \, \langle \uc_{(5,3)}, \uc_{(4,1)} \rangle 
- \, \frac{5}{13824} \, \langle \uc_{(5,3)}, \uc_{(3,2)} \rangle - 
\frac{1}{4320} \, \langle \uc_{(4,1)}, \uc_{(3,2)} \rangle,
\end{aligned} \label{C5.C11} \end{equation}
shows that $\uc_{(1,1)}$ can be recovered from the rest of the combinants. 

\subsection{} For the record, we state the symbolic expressions which were used to define the 
maps above. In formula~(\ref{C5.C32}), they are respectively 
\[ \begin{aligned} 
\Sc_{(5)} \otimes \Sc_{(5)} & \rightsquigarrow (a \,c \, \uu)^2 \, a_\ux^3 \, c_\ux^3, \\ 
\Sc_{(5)} \otimes \Sc_{(5,3)} & \rightsquigarrow (a \,c \,d) \, (a \,d \, \uu)^2 \, a_\ux^2 \, c_\ux^4, \\ 
\Sc_{(5)} \otimes \Sc_{(4,1)} & \rightsquigarrow (a \,c \, \uu) \, (a \,d \, \uu) \, a_\ux^3 \, c_\ux^3, \\ 
\Sc_{(5,3)} \otimes \Sc_{(5,3)} & \rightsquigarrow (a \,b \,d)^2 \, (a \,b \, \uu) \, (a \,d \, \uu) \, 
a_\ux \, c_\ux^5, \\ 
\Sc_{(5,3)} \otimes \Sc_{(4,1)} & \rightsquigarrow (a \,b \,d) \, (a \,b \, \uu)^2 \, a_\ux^2 \, c_\ux^4, 
\end{aligned} \] 
where the target of each map is $\Sc_{(8,2)}$. In (\ref{C5.C11}), they are respectively 

\[ 
\begin{aligned} 
\Sc_{(5)} \otimes \Sc_{(5,3)} & \rightsquigarrow (a \,c \,d)^2 \, (a \,d \, \uu) \, a_\ux^2 \, c_\ux^3, \\ 
\Sc_{(5)} \otimes \Sc_{(4,1)} & \rightsquigarrow  (a \,c \,d) \, (a \,c \, \uu) \, a_\ux^3 \, c_\ux^2, \\
\Sc_{(5)} \otimes \Sc_{(3,2)} & \rightsquigarrow  (a \,c \,d) (a \,d \, \uu) \, a_\ux^3 \, c_\ux^2, \\ 
\Sc_{(5,3)} \otimes \Sc_{(4,1)} & \rightsquigarrow (a \,b \,c) \, (a \,b \,d) \, (a \,b \, \uu) \, a_\ux^2 \, c_\ux^3, \\ 
\Sc_{(5,3)} \otimes \Sc_{(3,2)} & \rightsquigarrow (a \,b \,d)^2 \, (a \,b \, \uu) \, a_\ux^2 \, c_\ux^3, \\ 
\Sc_{(4,1)} \otimes \Sc_{(3,2)} & \rightsquigarrow (a \,b \,d) \, (a \,d \, \uu) \, a_\ux^2 \, c_\ux^3, 
\end{aligned} \] 
with target $\Sc_{(6,1)}$. 

\medskip 

{{\sc Acknowledgements:} 
\small 
The second author was partly funded by a discovery grant from NSERC. 
We are thankful to John Stembridge (author of the `SF' package for Maple). 
The G{\"o}ttinger Digitalisierungszentrum ({\bf GDZ}), 
the University of Michigan Historical Library ({\bf MiH}) as well as 
Project Gutenberg ({\bf PG}) have been useful in accessing some 
classical references.}

\bigskip
\centerline{---}

\vspace{1cm} 

\parbox{7cm}{\small 
{\sc Abdelmalek Abdesselam} \\ 
Kerchof Hall \\
Department of Mathematics\\
University of Virginia\\
P. O. Box 400137 \\
Charlottesville, VA 22904-4137\\ 
U.S.A. \\ 
{\tt malek@virginia.edu}}
\hfill \parbox{6cm}{\small 
{\sc Jaydeep Chipalkatti} \\ 
433 Machray Hall \\ 
Department of Mathematics\\
University of Manitoba \\ 
Winnipeg MB R3T 2N2 \\ Canada. \\ 
{\tt chipalka@cc.umanitoba.ca}}

\end{document}